\documentclass[16pt]{amsart}
\usepackage{graphicx} % Required for inserting images
\usepackage{amsmath}
\usepackage{amssymb,amsmath,amsthm,amsfonts,fullpage,float,latexsym, microtype, caption}
\usepackage{tikz}
\usepackage{pict2e}
\usepackage{centernot}
\usetikzlibrary{decorations.pathreplacing}
\usepackage{tikz-cd}
\makeatletter
\newtheorem*{rep@theorem}{\rep@title}\newcommand{\newreptheorem}[2]{%
\newenvironment{rep#1}[1]{%
\def\rep@title{\bf #2 \ref{##1}}%
\begin{rep@theorem}}%
{\end{rep@theorem}}}
\makeatother
\newreptheorem{theorem}{Theorem}

\newtheorem*{rep@proposition}{\rep@title}\newcommand{\newrepproposition}[2]{%
\newenvironment{rep#1}[1]{%
\def\rep@title{\bf #2 \ref{##1}}%
\begin{rep@proposition}}%
{\end{rep@proposition}}}
\makeatother
\newreptheorem{proposition}{Proposition}

\newtheorem{theorem}{Theorem}[section]
\newtheorem{proposition}[theorem]{Proposition}

\newtheorem{lemma}[theorem]{Lemma}
\newtheorem{corollary}[theorem]{Corollary}
\theoremstyle{definition}

\newtheorem{definition}[theorem]{Definition}
\newtheorem{observation}[theorem]{Observation}
\newtheorem{example}[theorem]{Example}

\newcommand{\xmathpalette}[2]{%
  \mathchoice
    {#1\displaystyle\textfont{#2}}
    {#1\textstyle\textfont{#2}}
    {#1\scriptstyle\scriptfont{#2}}
    {#1\scriptscriptstyle\scriptscriptfont{#2}}
}

\makeatletter

\newcommand{\mres@thickness}[1]{\dimexpr1.5\fontdimen8 #13\relax}

\newcommand{\mcorner}{\mspace{3mu}{\xmathpalette\mcorner@\relax}\mspace{3mu}}
\newcommand{\mcorner@}[3]{%
  \begingroup
  \setlength\unitlength{%
    \dimexpr0.8\fontcharht#21`A-0.5\mres@thickness{#2}% 80% height of capital letters
  }%
  \raisebox{0.5\dimexpr\mres@thickness{#2}}{%
    \begin{picture}(1,1)
      \roundcap\roundjoin
      \linethickness{\mres@thickness{#2}}% default rule thickness in the extension font
      \polyline(0,1)(0,0)(1,0)
    \end{picture}%
   }%
  \endgroup
}
\makeatother

\title{A Classification of Free and Free-Like Nilpotent Groups}
\author{Adam Moubarak}
\date{\today}

\begin{document}

\begin{abstract}
 Suppose $G$ is a $\mathcal{T}$-group (finitely generated torsion-free nilpotent) with centralizers outside of the derived subgroup being abelian of rank equal to $\text{rank}(Z_1)+1$. This includes the class of free nilpotent groups $\mathcal{N}_{r,c}$ of a given rank $r$ and class $c$. We show that the central series coincide in such groups and from this that they are metabelian. We then prove that all such groups arise as semidirect products of free abelian groups with respect to representation $[G,G]\to \text{UT}(n,\mathbb{Z})$ by automorphisms constructed from powers of elements in defining relations we call integral weights of $G$.
\end{abstract}
\maketitle

\newpage
\section{Introduction}
{\Large
Recall that a normal series of a group $G$
\begin{align}
    G=G_c\geq G_{c-1}\geq \dots \geq G_1\geq G_0=1,
\end{align}
is called central if $G_{i+1}/G_i\leq Z(G/G_i)$ - that is contained in the center of $G/G_i$ - for all $i$. Equivalently, for all $i$, $[G,G_{i+1}]\leq G_{i}$. A group having any central series is called nilpotent. Sometimes in the definition, a central series is not required to terminate, and in this case $G$ is nilpotent if any of its central series does. 
Two canonical central series are associated to any group: 
\begin{definition}
Define a series of $G$ inductively as $\gamma_1(G)=G$ and $\gamma_{i+1}(G)=[\gamma_i(G),G]$ for any $i\geq 1$. This is called the \emph{lower central series} of $G$. Requiring that the terms in the chain satisfy the latter condition for centrality of a series with equality.

Next, define the series $Z_0=1$, $Z_1=Z(G)$, and $Z_{i+1}/Z_i=Z(G/G_{i})$. This is called the \emph{upper central series} of $G$. Now, the first condition is satisfied with equality.
One can readily check that if (1) is any central series of $G$, then $\gamma_i\leq G_{c-i+1}\leq Z_{c-i+1}$ for each $i=1,\dots,c$.  In particular, if the LCS and UCS coincide, then there is only one central series of $G$.
\end{definition}
 It follows that $G$ is nilpotent if and only if its lower central or upper central series terminates in finitely many steps. Next, we gather common facts concerning the upper and lower central series of groups and those specific to nilpotent groups that are used freely and often cited.
Here, we use the following notation to express repeated commutators in $G$. Denote a repeated left-weighted $n$-commutator $$[\dots[[x_1,x_2],x_3],\dots,x_{n}]$$ by $[x_1,\dots,x_n]$, $x_i\in G$. At times, we will use right weighted commutators defined similary.

A finitely generated torsion-free nilpotent group $G$ is referred to as a $\mathcal{T}$-group. Notice that the center of $G$ is free abelian.
\begin{definition}
    Let $G$ be a $\mathcal{T}$-group. Suppose that the centralizer of every element $x\in G\setminus [G,G]$ is free abelian of rank $\text{rk}(Z_1)+1$. Then we say that $G$ has FL-centralizers or is an FL-centralizer group.
    
\end{definition}
Here FL stands for free-like motivated by the fact that free nilpotent groups of finite rank satisfy this property. In other contexts, this form of centralizers is also referred to as small in the literature.
It can be shown that this is equivalent to having for all $x\notin \gamma_2$ an element $u\notin \gamma_2$ with $C_G(x)=\langle u\rangle \times Z_1$ (see Appendix, Lemma \ref{FL-CentralizerCriterion}).
\begin{proposition}
Let $G$ be any group and $i$ and $j$ be positive integers. Then
\begin{enumerate}
    \item $[\gamma_i G,\gamma_j G]\leq \gamma_{i+j} G$.\\
    \item $[\gamma_i G,Z_j G]\leq Z_{j-i} G$ if $j\geq i$.
\end{enumerate}
\begin{proof}
(i). Fix $i$ and proceed by induction on $j$, with $j=1$ as base case holding trivially. Suppose $[\gamma_i,\gamma_j]\leq \gamma_{i+j}$. Note that $[\gamma_i,\gamma_{j+1}]=[\gamma_{j+1}, \gamma_i]=[\gamma_i, G, \gamma_i]$. Note that
\begin{align*}
    [G,\gamma_i,\gamma_j][\gamma_i,\gamma_j,G]\unlhd G
\end{align*}
and by the Three Subgroup Lemma $[\gamma_i,\gamma_{j+1}]$ is contained in this product. The induction hypothesis applied twice shows this product is contained in $\gamma_{i+j+1}$. Part (ii) follows similarly by induction on $i$.
\end{proof}
\end{proposition}
\begin{proposition}
If the center of a group $G$ is torsion-free, then each upper central factor is torsion-free. In particular, $G$ is torsion-free.
\end{proposition}
\begin{proof}
It suffices to show $Z_2/Z_1$ is torsion-free. Let $x\in Z_2$ and suppose $(xZ_1)^m=Z_1$ so that $x^m\in Z_1$ for some integer $m\neq 0$. For all $g\in G$, using an elementary identity and $[x,g]\in Z_1$, then $1=[x^m,g]=[x,g]^{x^{m-1}}[x^{m-1},g]=[x,g][x^{m-1},g]=\dots=[x,g]^m$. By assumption of torsion-freeness, this implies $x$ commutes with every element of $G$, so that $xZ_1$ is trivial in $G/Z_1$.
\end{proof}

\section{The Classification}

In this chapter, we determine the structure of $\mathcal{T}$-groups $G$ with FL-centralizers. We begin by proving that quotients by higher central terms are still FL-centralizer $\mathcal{T}$-groups and use this to give a short proof of coinciding central series for such $G$. Then, these results are applied to obtain several other corollaries eventually showing that $G$ is metabelian, as in the free nilpotent case. Finally, we construct explicit representations of $G$ into unitriangular groups $UT(n,\mathbb{Z})$ determined by collections of $r(c-3)$ independently chosen integers we call \emph{integral weights} of $G$. These gives a classification up to isomorphism (choices of embeddings). In examining explicit examples, the geometry of FL-centralizer nilpotent groups, more precisely their representation as groups of unitriangular matrices, becomes evident. 

\begin{proposition}
Let $G$ be a $\mathcal{T}$-group with FL-centralizers of class $c\geq 2$. Then the quotient $G/Z_{c-j+1}$ is also a $\mathcal{T}$-group with FL-centralizers of class $j-1$ for all $j=1,\dots,c$. 
\end{proposition}
\begin{proof}
Since the commutator subgroup of $G/Z_1$ is $\gamma_2(G)/Z_1(G)$ (no central element can have an FL-centralizer), it suffices to show that for all $g\in G$ and $h\notin \gamma_2$ such that $gZ_1$ commutes with $hZ_1$, we have that $g$ actually commutes with $h$ in $G$. Let $C_G(h)=\langle x\rangle \oplus Z_1$ and $h=x^iz$ for some nonzero integer $i$. Then 
\begin{align}
    [g,h]=[g,x^i]=[g,hx^{-1}][g,x]^{x^{i-1}}.
\end{align}
where $C_G(x)=C_G(h)$. We claim that $gZ_1$ also commutes with $xZ_1$ and that this makes the commutator vanish. %Consider the set of elements in $x\in G$ that solve an the equation of the form $x^i=k$ for some $k\in G$ (ie. the $i$-th roots). If $x_1,\dots,x_r$ are such elements that pairwise commute in $G$ then the subgroup generated by them is closed under taking $i$-th roots. It follows from considering $\langle x,z\rangle$ and divisibility of $Z_1(G)$ that 
%\begin{align*}
%    [g,x]&=[g,(hz^{-1})^{\frac{1}{i}}]\\
%    &=[g,ab]\\
%    &=[g,a],
%\end{align*}
%for unique elements $a=h^{\frac{1}{i}}$ and $b=(z^{-1})^{\frac{1}{i}}$. First consider the case $c=2$. Note $[g,h^{\frac{1}{i}}]^{1\cdot i}=[g,h]$ which implies 
%\begin{align*}
 %   [g,x]=[g,h^{\frac{1}{i}}]\in Z_1.
%\end{align*}
%so that $(h^{\frac{1}{i}})^{-1}hg^{-1}=g^{-1}(h^{\frac{1}{i}})^{-1}h$ or equivalently $h^{i-1}(g^{-1})^i=(g^{-1})^{i}h^{i-1}$. For $c>2$, the proof is similar choosing any central elements $z_1,\dots,z_{c-2}$ of $G$ then $[g,x]^{i^{c-1}}=[g,h^{\frac{1}{i}},z_1^{\frac{1}{i}},\dots,z_{c-2}^{\frac{1}{i}}]^{i^{c-1}}\in Z_1$. Then equation (1) reads
%\begin{align}
 %   [g,h]=[g,x^{-1}][g,x]
%\end{align}
%As mentioned above, this argument holds with the replacement $h\to h^{-1}$ and $x\to x^{-1}$, which by symmetry of the equation means $[g,h]=[g,h^{-1}]$ or that $g\in C_G(h)$.

Notice that (1) can be rewritten as follows: first
\begin{align*}
    [g,hz^{-1}x^{-1}]=[g,hx^{-1}]=[g,x^{-1}][g,h]^{x^{-1}}
\end{align*}
so that
\begin{align*}
    [g,h]&=[g,h]^{x^{-1}}[g,x^{-1}][g,x]^{x^{i-1}}\\
    \implies 1&=[g,x^{-1}][g,x]^{x^{i-1}}
\end{align*}
by centrality of $[g,h]$. At this point use the symmetry $x\to x^{-1}$, $i\to -i$ in the equation for $h$ to get
\begin{align*}
    1&=[g,x][g,x^{-1}]^{x^{i+1}}\\
    &=[g,x]\left([x,g]^{x^{i-1}}\right)^{x^{i+1}}\\
    &=[g,x][x,g]^{x^{2i}}.
\end{align*}
This means $[x,g]$ commutes with $x^{2i}$, hence with $x$. Letting $gx^{-1}g^{-1}=x^rc$ for some nonzero integer $r$ and central $c \in G$ we have $[x,g]=x^{r+1}c$ and $[g,x^{-1}]=x^{r+1}c$ so that equation (1) reads $[g,h]=1$. The rest follows by induction on $j$.
\end{proof}
\begin{corollary}
Let $G$ be a $\mathcal{T}$-group with FL-centralizers of class $c\geq 2$. Then the central series of $G$ coincide.
\end{corollary}
\begin{proof}
Proceed by induction on $c$ with inductive step as follows. Suppose all $\mathcal{T}$-groups with FL-centralizers of class less than $c$ have coinciding series and $G$ have class $c$. The group $G/Z_1$ has coinciding central series by the proposition and hypothesis. That is, we have the identity $$Z_{c-i}(G/Z_1)=\gamma_i(G/Z_1)$$ for all $i=1,2,\dots,c-1$. The left hand side is equal to $Z_{c-i+1}/Z_1$ and right hand side is a well-defined group $\gamma_i/Z_1$ using that $Z_1=\gamma_c\leq \gamma_i$ for any group with FL-centralizers. Whence, for each $i=1,\dots,c$, we have
\begin{align*}
    Z_{c-i+1}\leq \gamma_iZ_1\leq \gamma_i.
\end{align*}
\end{proof}

An alternate proof that the upper and lower central series of $G$ coincide is provided in the appendix.
In general a nilpotent group of class $c$ can have central series of length greater than $c$ as well as many distinct central series of a given length. If the lower and upper central series of $G$ coincide then it has only one central series of length $\leq c$ thereby simplifying the study of the central series the group admits.

It is well-known that all finite $p$-groups and all free nilpotent groups of any finite rank and class have equivalent upper and lower central series. It is also known that the non-free nilpotent unitriangular groups $\text{UT}(n,\mathbb{Z})$ also have coinciding series. Given an $m$-generated free $c$-nilpotent group $R=F/\gamma_{c+1}(F)$, where $F$ is an absolutely free group on $m$ generators, we must have $Z_{c-k+1}(R)=\gamma_k(R)$ for each $k=1,\dots, c$. We prove the same equivalence for the class of nilpotent groups consisting of those with FL-centralizers.
 However, there is still no hope that this will exhaust all $\mathcal{T}$-groups with coinciding series as the following shows:
\begin{observation}
For any integer $k\geq 4$, there exists a proper (non-free) subgroup $H\leq \text{UT}(k,\mathbb{Z})$ that has coinciding central series but not $FL$-centralizers.
\end{observation}
\begin{proof}
Given $k$, let $\tilde{H}=\text{UT}(n,\mathbb{Z})$ where $3\leq n< k$ which has equal Upper and Lower Central Series but it is readily seen not to have FL-centralizers. There is a natural embedding $J:\text{UT}(n,\mathbb{Z})\to \text{UT}(k,\mathbb{Z})$ filling in the top (k-n+1) - superdiagonals, then the image subgroup of $\tilde{H}$ under $J$ is such an example.
\end{proof}
We borrow familiar definitions from the literature. A group $G$ is called an \emph{$R$-group} if whenever roots of elements exist, they are unique: that is, for all $a\in G$ and integers $n\in\mathbb{Z}$ the equation $x^n=a$ has at most one solution in $G$. Examples include all torsion-free nilpotent groups, in particular $\mathcal{T}$-groups.

A subgroup $H$ of an $R$-group is said to be \emph{isolated} if for every $a\in H$, all roots of $a$ in $G$, provided they exist, also lie in $H$. The isolated closure of a subset $S$ of an $R$-group $G$ is defined as the minimal isolated subgroup of $G$ containing $S$. In this terminology, a subgroup is isolated if it coincides with its isolated closure. Notice an isolated subgroup $H=\langle x_1,\dots,x_n\rangle \leq G$ may not have the property that each generator has a root, since one may not exist in $G$.
We next consider many special cases. Suppose $x_1,\dots,x_n\in G$ commute with each other and all have roots of a common order, say $k\in\mathbb{N}$. Denote by $\alpha_i=x_i^{\frac{1}{k}}$ the unique element of $G$ with $\alpha_i^k=x_i$. Then the abelian subgroup $\langle x_1,\dots,x_n\rangle$ of $G$ is such that each of its elements has a root, of order $k$, in $G$. In fact the roots themselves form a subgroup $\langle \alpha_1,\dots,\alpha_n\rangle < G$ which is contained in $\langle x_1,\dots,x_n\rangle$ when the latter is isolated. However, in this case, we have $x_j=x_{1}^{r_1}\dots x_{n}^{r_n}$ for some integers $r_i(j)$ with common divisor $k>1$. If the $\{x_i\}$ form a free basis for $H$ then none could have roots of the same order in $G$. The reason for this is that the above relation would require $r_j=k=1$. Hence, if each basis element has a root, generation is bounded by restrictions on the exponents. 

Since $G$ is torsion-free and by assumption the generating set is minimal, so no $x_i$ is a root of $x_j$ for $i\neq j$, we must have at least two nonzero exponents (if $r_1=1$ we must have at least \emph{three} nonzero exponents).

Even when there are no common order roots of the collection of elements, we can still salvage an isolated subgroup of $\langle x_1,\dots,x_n\rangle$: if $x_i$ has a $k_i$-th root that satisfies $\alpha_i^{k_i}=x_i$ in $G$ for each $i$ and $m=k_1\dots k_n$ then each power $x_i^{\frac{m}{k_i}}$ has an $m$-th root given by $\alpha_i^{m}$ so that an isolated subgroup containing $x_1,\dots,x_n$ is the power subgroup $\langle x_1^m,\dots,x_n^m\rangle$ (not necessarily their isolated closure, however). The above result gives that all subgroups $\gamma_j=Z_{c-j+1}$ are isolated in $G$ and much of the machinery in the proofs depends on the commutator subgroup being isolated in particular.
\subsection{The Structure of the Commutator Subgroup}
In this section, we study properties of $\gamma_2$ for a $\mathcal{T}$-group $G$ with FL-centralizers. To do this, we make use of the canonical homomorphisms usually constructed in Grun's Lemma on perfect groups. More generally, these types of maps appear in tensor product constructions for nilpotent groups. 

We informally use the term separability for any property of a group that involves trivial intersection of subgroups related by some automorphism of $G$. Examples of properties that fall under this are malnormality of subgroups and more generally of being conjugacy stable abelian (CSA). Another type of separability conditions considered are those involving the preservation of "separated elements" of a group under homomorphism surjecting to finite groups. An example of this is the Conjugacy-separable property of a group in which separated here means lying in distinct conjugacy classes. It is well known that the property of being CSA is stronger than CT (commutative transitivity) for any group. It can be shown that any $\mathcal{T}$-FL-centralizer group $G$ is \emph{locally} CT by means of a defined divisibility relation on $G$ and this is analogous to the fact having abelian centralizers is equivalent to CT on the whole group (see the Appendix for more details).

What we are really saying is a non-abelian FL-centralizer group $G$ has commutative transitivity if commutativity is transitive on the complement $G-\gamma_2$. This route leads to showing that the maximal abelian subgroups of $G$ are centralizers $C_{G}(x)$ for $x\notin \gamma_2$ and the commutator subgroup of $G$.
\begin{proposition}
Let $G$ be a $\mathcal{T}$-group with FL-centralizers. Then the only maximal subgroups of $G$ are centralizers outside the derived subgroup and these are malnormal. 
\end{proposition}
\begin{proof}
Given any $h\in G\setminus \gamma_2$ let $C_G(h)=\langle u\rangle \times Z_1$. We claim $g^{-1}C_G(h)g\cap C_G(h)=1$ for all $g\in G\setminus{C_G(h)}$. Indeed if $gu^{k_1}z_1g^{-1}=u^{k_2}z_2$ where $k_1,k_2$ are integers and $z_1,z_2$ central elements then
\begin{align*}
    gu^{k_1}g^{-1}=u^{k_2} \: \text{mod} \: Z_1\\
    [g,u^{k_1}]=u^{k_2-k_1} \: \text{mod} \: Z_1
\end{align*}
As $Z_{c-1}$ is isolated and $u^{k_2-k_1}\in Z_{c-1}$ we must have $k_1=k_2$. But then $G/Z_1$ is an $R$-group as a torsion-free nilpotent group so $gu^{k_1}g^{-1}=u^{k_1}$ implies $g\in C_G(u)=C_G(h)$.

Next we show the centralizers are maximal abelian. Suppose $A$ is any maximal abelian subgroup with $C_G(x)\subseteq A$. In particular $A$ cannot be contained in the derived subgroup of $G$ so let $h\in A-\gamma_2$ and $C_G(h)=\langle v\rangle \times Z_1$. Since $x\in A$ and the latter is commutative we have $C_G(h)=C_G(x)$. That is, we have
\begin{align*}
    A\gamma_2/\gamma_2\leq C_G(x)\gamma_2/\gamma_2
\end{align*}
or, by taking the preimage of both sides under the map $G\to G^{\text{ab}}$, we have $A\gamma_2\leq C_G(x)\gamma_2$. First notice that $C_G(x)\cap \gamma_2=Z_1$ otherwise if $\omega\in \gamma_2$  commutes with $x$ then $\omega=v^kz$ which implies $k=0$, as $v\notin \gamma_2$, and so $\omega=z\in Z_1$. Thus we have an isomorphism
\begin{align*}
    A\gamma_2/\gamma_2\subseteq C_G(x)\gamma_2/\gamma_2\approx C_G(x)/\gamma_2\cap C_G(x)=C_G(x)/Z_1=\langle v\rangle
\end{align*}
so that $A\gamma_2/\gamma_2$ is cyclic and tracking through the mappings we have $A\gamma_2/\gamma_2=\langle v^m\gamma_2\rangle$ for some integer $m$. Hence, for all $a\in A$, there is an integer $m'$ and $c\in \gamma_2$ such that $a=v^{m'}c$. Then since $A$ is abelian $c\in Z_1$ and we conclude $a\in C_G(x)$.

To see why these are the only maximal abelian subgroups of $G$, if $A$ is such a subgroup and there is $x\in A\setminus \gamma_2$ then $C_G(x)\leq A$ and the above shows $A=C_G(x)$. On the other hand, if $A\leq \gamma_2$ then 
\begin{align*}
    A=C_G(A)=\bigcap_{a\in \gamma_2} C_G(a)=\gamma_2.
\end{align*}
\end{proof}
This is saying that, at least when the centralizers are small, the property of being CSA is "local" in $G$. It is natural to ask whether all maximal abelian subgroups of $G$ are centralizers of elements outside the derived subgroup or in this sense "global". In fact, it would follow from the previous result that $G$ is CSA. However, we show this cannot be the case by proving the stronger condition of being metabelian.
\begin{lemma}\label{refinement of Grun's Lemma}
For a $\mathcal{T}$-group $G$ with FL-centralizers of class $c\geq 3$ we have
\begin{align*}
    G^{\text{ab}}\simeq \langle [a,G] \rangle
\end{align*}
for every $a\in Z_2\setminus Z_1$. From the isomorphisms this is equivalent to
\begin{align*}
    C_G(Z_2)=\gamma_2.
\end{align*}
\end{lemma}

In particular, the homomorphisms obtained from Grun's Lemma all produce isomorphic subgroups of $Z_1(G)$. 
\begin{proof}
For each $a\in Z_2\setminus Z_1$ consider the nontrivial homomorphism $\varphi_a:G\to Z_1(G)$ given by $[a,x]$. In any group we have $\gamma_2\subset \text{ker}(\varphi_a)$ as $[\gamma_2,Z_2]\leq Z_{2-2}=\{1\}$. However, equality fails in general even for nilpotent groups. For any $x\notin \gamma_2$ if $[a,x]=1$ then $a\in \langle u\rangle \times Z_1$ so that $a=u^iz$ for some $u\notin \gamma_2$ and $z\in Z_1$ which gives $u^i\in Z_2=\gamma_{c-1}\subset \gamma_2$. It must then be that $a\in Z_1$, a contradiction. Hence no element outside the derived subgroup of $G$ can commute with $a$ which means $\text{ker}(\varphi_a)=\gamma_2$.

Equivalence with the other statement follows immediately:
\begin{align*}
    C_G(Z_2)=\left(\bigcap_{a\in Z_2\setminus Z_1}C_G(a)\right)\cap G=\gamma_2.
\end{align*}
The remaining statement is also clear as for any two non-central elements $a,a'\in Z_2$, letting $\varphi_x^*$ denote the isomorphism obtained by passing to the abelianization of $G$ we get
\begin{align*}
    \varphi_a^*(G^{ab})\simeq G^{ab} \simeq \varphi_{a'}^*(G^{ab}).
\end{align*}
\end{proof}
Actually, the central subgroup $\langle [a,G]\rangle$ is written for simplicity of notation but actually the stronger statement is true that the single commutators $\{[a,x]\}_{x\in G}$ generate the group. Given a basis of the quotient $G^{ab}$, pulling back along the abelianization and pushing forward along the homomorphism $x\mapsto [a,x]$ gives a finite set of the form $\{[a,x_1]\dots,[a,x_r]\}$.

It follows from the previous result that
\begin{theorem}
    Any $\mathcal{T}$-group $G$ with FL-centralizers is metabelian.
\end{theorem}
\begin{proof}
If $G$ is 3-step nilpotent the above facts show $\gamma_2=C_G(Z_2)=C_G(\gamma_2)$ so that $G$ is metabelian. Suppose that all FL-centralizer $\mathcal{T}$-groups of class $c-1$ are metabelian and let $G$ be such a group of class $c$. %First notice that application of Lemma \ref{refinement of Grun's Lemma} to $\tilde{G_k}=G/Z_{c-k+1}$ gives the following: for any $a\in Z_{c-k+1}\setminus Z_{c-k}$ with arbitrary $k\geq 1$, since $Z_{c-k+1}/Z_{c-k-1}=Z_2(G/Z_{c-k-1})$ and $Z_{c-k}/Z_{c-k-1}=Z_1(G/Z_{c-k-1})$, so the element $aZ_{c-k-1}$ in $G/Z_{c-k-1}$ has centralizer 
%\begin{align*}
 %   C_{G/Z_{c-k-1}}(aZ_{c-k-1})=\gamma_2(G)/Z_{c-k-1}.
%\end{align*}
Fix any $a\in \gamma_2$ and let $y \in \gamma_2$ arbitrary. Then, the group $G/Z_1$ is metabelian so that its commutator subgroup $\gamma_2/Z_1$ is abelian. In particular, we have $[a,y]\in Z_1$ and
\begin{align*}
    a[a,y]&=[a,y]a\\
    y^{-1}ay&=a^{-1}(y^{-1}ay)a
\end{align*}
or in other words $a^y\in C_G(a)$. Moreover, the element $[a,y]$ being central, we can write $a^y=az$ for some $z\in Z_1$. It follows that $C_G(a^y)=C_G(a)$. The following computation can then be made
\begin{align*}
    [a,y]&=[a^y,y]=(a^y)^{-1}y^{-1}a^yy\\
    &=y^{-1}a^{-1}y^{-1}ay^2=[y,a]y^2\\
    [a,y]^2&=y^2
\end{align*}
which by unique root extraction in $G$ means $y\in Z_1$ and in particular $[y,a]=1$. Thus, for all $\omega_1,\omega_2 \in \gamma_2$ we have
\begin{align*}
    C_G(\omega_1)=\gamma_2=C_G(\omega_2).
\end{align*}
\end{proof}

\begin{corollary} (Logarithim operation)
    For any $a\in \gamma_2$ we have
\begin{align*}
        a^x=a^y\implies x=y \: (\text{mod}\: \gamma_2)
\end{align*}
\end{corollary}
\begin{proof}
    It sufficies by induction to consider only $a\in Z_2\setminus Z_1$. Let $x,y\in G- \gamma_2$ and $C_G(x)=\langle u\rangle\times Z_1$ and $C_G(y)=\langle v\rangle \times Z_2$. Put $x=u^iz$ and $y=v^jz_1$ for integers $i,j$ and $z,z_1\in Z_1$. Notice $[a,x]=[a,y]$ so that
    \begin{align*}
        [a,u^i]=[a,v^j]
    \end{align*}
  which implies by the previous result that there is $c\in \gamma_2$ with $v^j=cu^i$. Then
  \begin{align*}
      y=cu^iz_1=cxz^{-1}z_1
  \end{align*}
  and we take this equation modulo $\gamma_2$.
\end{proof}

\subsection{Integral Weights of a Group}
\begin{theorem}\label{Weighted Roots Classification}
Let $G$ be a $\mathcal{T}$-group with FL-centralizers of class $c>3$. Then the structure of $G$ is determined as follows: there exist a collection of $(c-3)$ integers $(k_{34},\dots,k_{(c-1,c)})$ called the weights of $G$ and a collection of elements $\{u_3,\dots,u_c\}\subset G\setminus [G,G]$ such that $z_{ij}\equiv u_ju_i^{-k_{ij}}\in Z_{c-i+2}$ for all $j>i$, that satisfy the following relations: for all $3\leq i<j<t\leq c$, we have
\begin{enumerate}
    \item Multiplicative Pascal-relation: $k_{ij}k_{jt}=k_{it}$\\
    \item Multiplication in $Z_{c-i+2}$ given by: $z_{ij}^{k_{jt}}z_{jt}=z_{it}$
\end{enumerate}
\end{theorem}
The collection $(k_{34},\dots,k_{c-1,c})$ called the \emph{integral weights} of $G$ and $(c-3)$ is the \emph{weighted dimension} of $G$.

\begin{proof}
Fix $h\in G\setminus \gamma_2$ and set $\bar{h_j}=hZ_{c-j+1}$ for all $j=1,\dots,c$. Since each $\bar{G_j}=G/Z_{c-j+1}$ has FL-centralizers there exist elements $\bar{u_j}\in \bar{G_j}\setminus \gamma_2(\bar{G_{j}})$ (so that $u_j\notin \gamma_2$) such that
\begin{align*}
    C_{\bar{G_j}}(\bar{h_j})=\langle \bar{u_j}\rangle \times \bar{Z_{c-j+2}}.
\end{align*}
In words $C_{\bar{G_j}}$ consists of all cosets with representative in $G$ whose commutator with $h$ falls in $Z_{c-j+1}$. It is thus clear that representatives of elements in $C_{\bar{G_j}}$ also appear as representatives of elements in $C_{G_i}$ for all $i<j$.

If $\Omega_k:G\to G/Z_{c-k+1}$ is the natural projection for all $k$, this condition is expressed as $
    \Omega_j^{-1}(C_{\bar{G_j}})\leq \Omega_i^{-1}(C_{\bar{G_i}})
$ for all $i<j$. That is, we obtain a chain of subgroups of $G$:
\begin{align*}
    \mathcal{C}_h: 1\leq \Omega_c^{-1}(C_{\bar{G_c}})\leq \dots \leq \Omega_3^{-1}(C_{\bar{G_3}})\leq G
\end{align*}
for each $h\in G\setminus \gamma_2$. Note that $\Omega_j^{-1}(C_{\bar{G_j}})=\langle u_j, Z_{c-j+2}\rangle$ for all $j=1,\dots, c$ and the series is allowed (despite length $<c$) because the terms are not normal in $G$; otherwise $u_j$ would lie in some proper higher center of $G$. This also gives the relations $u_j=u_i^{k_{ij}}z_{ij}$ for an integer $k_{ij}\neq 0$ and elements $z_{ij} \in Z_{c-i+2}$ over all $j>i\geq 3$. Likewise for any $j$ consider all indices $j<t$ with associated relation 
\begin{align*}
    u_t&=u_j^{k_{jt}}z_{jt}=(u_i^{k_{ij}}z_{ij})^{k_{jt}}z_{jt}\\
    &=u_i^{k_{ij}k_{jt}}z_{ij}^{k_{jt}}z_{jt}
\end{align*} 
comparing with $u_t=u_i^{k_{it}}z_{it}$ and using that $\langle u_i\rangle\cap Z_{c-i+2}=1$ gives the desired relations.
\end{proof}

Notice this allows us to define a representation $L_H:G^{\text{ab}}\to \text{UT}(c-2,\mathbb{Z})$ with respect to a given free basis. Explicitly, select any basis $\{\bar{h_1},\dots,\bar{h_r}\}$ for the abelianization of $G$ where $\bar{h_i}=h_i\gamma_2$ for each $i=1,\dots,r=\text{rk}(G)$. For each $\ell=1\dots,r$, suppose $(k_{ij}^{\ell})_{i<j}$ are the weighted roots of $h_{\ell}\in G$. Then define the matrix 
\begin{align*}
    \left(L(h_{\ell})\right)_{(i-2,j-2)}=1+k_{ij}e_{i-2,j-2}
\end{align*} for all $i<j$ between $3$ and $c$. Now we can extend this map homomorphically to $G^{ab}$ since each element has a unique expression in this basis: $L(\alpha_1 \bar{h_1}+\dots+\alpha_r \bar{h_r})=L(h_{1})^{\alpha_1}\cdots L(h_{r})^{\alpha_r}$.
In general, the kernel is nontrivial and the Pascal-relation among the entries $k_{ij}k_{jk}=k_{il}$ of $L$ can be used to describe it. Actually, there is a rich underlying symmetry since in general each element $h^{\alpha_1}\dots h^{\alpha_r}\in \text{ker}(L)$ for integers $\alpha_i$ will give a system of polynomial equations in the $r(c-3)$ variables $(k_{34}^j,\dots,k_{(c-1,c)}^j)_{h_j}$, indexed over $j=1,\dots,r$.

Hence we realize the abelianization of each $\mathcal{T}$-group with FL-centralizers, in particular free nilpotent groups, as a subgroup of a unitriangular group $\text{UT}(c-2,\mathbb{Z})$. Precomposing with the abelianization map gives a representation
\begin{align*}
    \varphi:G\to \text{UT}(c-2,\mathbb{Z}).
\end{align*}
which depends on the choice of free basis, a priori, whose kernel is precisely $\gamma_2$.  However, notice that changing basis can be viewed as a sequence of Nielsen-transformations each of which is an elementary transformation of the associated matrix whose columns are given in the $\{h_1,\dots,h_r\}$ basis. The two image subgroups under $L$ and $L'$ are related by conjugation by the resulting element of $GL(r,\mathbb{Z})$. Here we make dimensions can be consistent by embedding in $\text{UT}(c-2,\mathbb{Z})$ if $r+2\leq c$ or vice versa for $c\leq r+2$. Hence, the choice of basis is immaterial up to isomorphism.

Recall for any group $G$ there is a natural homomorphism 
\begin{align*}
    \text{Aut}(G)\to \text{Aut}(G^{\text{ab}}).
\end{align*} It is known that if $F_r$ is a free group on $r$ generators then this map is surjects onto $\text{GL}(r,\mathbb{Z})$. For a free nilpotent or solvable group $G$ of rank $r$ and class $c$, by passing to the quotient by its relator subgroup, the map $\text{Aut}(G)\to \text{GL}(r,\mathbb{Z})$ turns out to be surjective as well.

Conversely, suppose we have a collection of $r$ sets of $(c-3)$ integers $(k_{ij})$ with each satisfying the relation and consider the subgroup $U\leq \text{UT}(c-2,\mathbb{Z})$ generated by the associated unitriangular matrices. We will have a natural isomorphism  $G^{\text{ab}}\approx \langle L(h_1),\dots,L(h_r)\rangle$. To reconstruct $G$ we find splitting extensions, using freedom in the choice of structure of its commutator subgroup.
\subsection{Examples}
In this section, we use the previous scheme to quite explicitly construct the $\mathcal{T}$-groups with FL-centralizers of any given rank $r$ and step $c\geq 3$ up to isomorphism.
\begin{theorem}\label{The Classification}
Up to isomorphism, every rank $r$ and step $c$ $\mathcal{T}$-group with FL-centralizers is given by
\begin{align*}
    G\approx \mathbb{Z}^{a}\rtimes_{\Phi} \langle A_1,\dots,A_r\rangle
\end{align*}
where $R=\langle A_1,\dots,A_r\rangle\leq \text{UT}(c-2,\mathbb{Z})$ has generators 
\begin{align*}
    A_i=\begin{pmatrix}
        1 & k_{34}^i & k_{34}^ik_{45}^i & \cdots & k_{34}^i\dots k_{c-1,c}^i\\
        0 & 1 & k_{45}^i & \cdots & k_{45}^i\dots k_{c-1,c}^i\\
        0 & 0 & 1 & \ddots & k_{c-1,c}^i\\
        0 & 0 & 0 & \cdots & 1
    \end{pmatrix}
\end{align*}
for some collections of integers $(k_{34}^i,\dots,k_{c-1,c}^i)$, $i=1,\dots r$.
\end{theorem}
\begin{proof}
Any such $G$ has a short-exact sequence
\begin{align*}
    1\to \mathbb{Z}^{a} \to G \to G^{ab}\to 1 
\end{align*}
where the inclusion is of the abelian commutator subgroup $[G,G]\approx \mathbb{Z}^a$ using that $G$ is metabelian. It suffices to produce a homomorphism 
\begin{align*}
    \Phi: G^{ab}\to \text{Aut}(F)
\end{align*}
into the automorphism group of some group $F$ isomorphic to $[G,G]$. A natural choice is to take $F=\mathbb{Z}^a$, for integer $a=\text{rank}([G,G])$. Of course any other choice of $F$ gives the same automorphism group up to isomorphism so our representations $\Phi$ of $G$ are effectively equivalent. The simplest choice is the the composition $$\Phi=\iota\circ L:G^{ab}\to GL(c-2,\mathbb{Z})$$ where $L: G^{ab}\to UT(c-2,\mathbb{Z})$ is the embedding given by weights of $G$ in Theorem \ref{Weighted Roots Classification} and the latter is $UT(c-2,\mathbb{Z})$. The converse also follows that such groups are FL-centralizer $\mathcal{T}$-groups.
\end{proof}
 Suppose one desires to find an explicit embedding of $G$ into a unitriangular group $UT(n,\mathbb{Z})$. We may first begin by choosing an embedding of the first factor $\mathbb{Z}^{a}$ of the semi-direct product, so that the underlying set is a cartesion product of subsets of matrices of the same dimension, as in general the free abelian ranks of $[G,G]$ and $G^{ab}$ may differ. To do this, we need at least $a$ dimensions of freedom in the superdiagonal elements which imposes a lb on dimension
 \begin{align*}
     a\leq \frac{n(n-1)}{2}
 \end{align*}
Let $\psi:\mathbb{Z}^a\to UT(n,\mathbb{Z})$ be an explicit embedding which identify with the isomorphism onto its image. We will show that all choices lead to isomorphic groups. Indeed,
\begin{proposition}
Let $A,B$ be any groups with $\Phi:B\to\text{Aut}(A)$ a homorphism given by $b\to \phi_{b}$, and $\Psi:A\to C$ be an isomorphism. Then the following semi-direct products are naturally isomorphic:
\begin{align*}
    A\rtimes_{\Phi}B\simeq C\rtimes_{\Psi\circ \Phi\circ \Psi^{-1}}B 
\end{align*}
\end{proposition}
Letting $C_{e_i}\in UT(n,\mathbb{Z})$ denote the image of free generator $e_i=(0,\dots,1,\dots,0)\in \mathbb{Z}^a$ for each $i=1,\dots,a$, and the abelian subgroup $C$ generated by them, we have as sets $G\subset UT(n,\mathbb{Z})\times UT(r,\mathbb{Z})$. The latter naturally sits in $UT(n+r+2,\mathbb{Z})$ but one must find a mapping that respects the semi-direct group operation on $G$. 

\begin{example}
There is a shortcut to constructing 3-step nilpotent groups with FL-centralizers in the sense that we do not need the metabelian result, suggested as follows: let $G$ be any FL-centralizer group of step $c>3$ and consider that Lemma \ref{refinement of Grun's Lemma} gives 
\begin{align*}
    C_G(Z_2)=\gamma_2, C_{G/Z_1}(Z_3/Z_1)=\gamma_2/Z_1,\dots,C_{G/Z_{c-3}}(Z_{c-1}/Z_{c-3})=\gamma_2/Z_{c-3}.
\end{align*}
In particular, the quotient $Z_{c-1}/Z_{c-3}$ is an abelian  subgroup of $G/Z_{c-3}$ so is torsion-free. Then it is natural ask whether the short exact sequence
\begin{align*}
1\to Z_{c-1}/Z_{c-3}\to G/Z_{c-3}\to G/Z_{c-1}\to 1
\end{align*}
splits. This occurs precisely when we can find a homomorphism $\Phi:\mathbb{Z}^q\to \text{GL}(p,\mathbb{Z})$. For instance, when $(p,q)=(3,2)$ send the basis vectors of $\mathbb{Z}^2$ to 
\begin{align*}
    \Phi(e_1)=\begin{pmatrix}
        1 & 1 & 0\\
        0 & 1 & 0\\
        0 & 0 & 1
    \end{pmatrix},
    \Phi(e_2)=\begin{pmatrix}
        1 & 0 & 0\\
        0 & 1 & 1\\
        0 & 0 & 1
    \end{pmatrix}
\end{align*}
Extending this to a homomorphism $\Phi:\mathbb{Z}^2\to GL(3,\mathbb{Z})$ gives that
\begin{align*}
    \Phi(ae_1+be_2)=\begin{pmatrix}
        1 & a & ab\\
        0 & 1 & b\\
        0 & 0 & 1
    \end{pmatrix}
\end{align*} and we see the multiplicative relations between integral weights appear in the products of elements in $G$! Set $G=\langle f_1,f_2,f_3\rangle \rtimes_{\Phi}\langle e_1,e_2\rangle$ where $f_i$ are the standard basis vectors of $\mathbb{Z}^3$. Then $G$ is an FL-centralizer nilpotent group with operation
\begin{align*}
    x\cdot y&=((a_1,a_2,a_3),(\varepsilon_1,\varepsilon_2))\cdot ((b_1,b_2,b_3),(\delta_1,\delta_2))\\
    &=((a_i)+\Phi_{(\varepsilon_i)}(b_i),(\varepsilon_i+\delta_i))
\end{align*}
where $\Phi_{(\varepsilon_i)}(b_i)=\Phi(\varepsilon_1,\varepsilon_2)(b_1f_1+b_2f_2+b_3f_3)=\begin{pmatrix}
    1 & \varepsilon_1 & \varepsilon_1\varepsilon_2\\
    0 & 1 & \varepsilon_2\\
    0 & 0 & 1
\end{pmatrix}\begin{pmatrix}
    b_1\\
    b_2\\
    b_3
\end{pmatrix}$ so that
\begin{align*}
    x\cdot y=\left( \begin{bmatrix}
        (a_1+b_1)+\varepsilon_1b_2+\varepsilon_1\varepsilon_2b_3\\
        (a_2+b_2)+\varepsilon_2b_3\\
        (a_3+b_3)
    \end{bmatrix},\begin{bmatrix}
        \varepsilon_1+\delta_1\\
        \varepsilon_2+\delta_2
    \end{bmatrix}\right)
\end{align*}
\end{example}

\subsection{Geometry of Linear Representations}
Geometrically, we see that in the first component, the 3-dimensional vector $\vec{a}$ is added to a transformed version of $\vec{b}$ according to the discrete weights of the group. The addition of z-components is unaffected, scaled in the y-direction by $\varepsilon_2$ (left or right depending on sign) and then in the x-direction scaled based on the scaling in $y$. The second component just adds the 2-dimensional vectors. In fact this interpretation of FL-centralizer groups is parallel to other geometric treatments of nilpotent groups, particularly free ones. The above example suggests that distortion in $G$ is measured by the vectors
\begin{align*}
    \begin{bmatrix}
        \varepsilon_1b_2+\varepsilon_1\varepsilon_2b_3\\
        \varepsilon_2b_3
    \end{bmatrix}
\end{align*}
Actually, it and the weights by extension measure deviation from commutativity in $G$ which the commutator subgroup does too and abelian groups have no distorted subgroups. In this sense, we the amount of distorted subgroups and how distorted one may expect to be tracked by the weights as well.

We also expect freeness of a nilpotent group to be be contained in the integral weights of $G$. From this previous section we see that the structure of $G$ is determined by its abelianization and its derived subgroup. Due to this fact, it is therefore of interest to intepret the commutator elements geometrically. By the classification theorem, we saw that each element of the abelianization $G^{ab}$ represents an automorphism of a discrete lattice $\mathbb{Z}^a$ where $a=\text{rk}([G,G])$. However, we also saw that $G^{ab}\simeq \langle [a,G]\rangle$ for each $a\in \gamma_2$ by a refinement of Grun's Lemma. Hence, the structure of the derived subgroup alone ultimately determines $G$. 
\begin{proposition}
A free nilpotent group $G$ of class $c$ and rank $r$ with the $R_{\infty}$ property satisfies
\begin{align*}
    2r-1\leq \text{rk}(\gamma_2).
\end{align*}
\end{proposition}
\begin{proof}
By Lemma 2.3 [Escayola, Rivas], we have $c\leq \text{rk}(\gamma_2)+1$ taking $A=\gamma_2$. This follows from the result that $\gamma_2$ is maximal abelian subgroup of $G$. It has been shown that when $G$ is free nilpotent it has the $R_{\infty}$ property precisely when $c\geq 2r$. Suppose we had $\text{rk}(\gamma_2)<2r-1$. Then, by the above inequality
\begin{align*}
    c< 2r-1+1=2r.
\end{align*}
\end{proof}

\newpage

\section{Appendix}

To motivate the results, it is useful to have criterion for a nilpotent group to have FL-centralizers. We specifically look for criterion expressed in terms of equation(s) in the group; it can be expressed in model theoretic language as a translation of the FL-centralizer property into a logical equivalence between certain first and second order sentences in the group. This follows an analogous construction of divisibility and associates- elements differing by a unit- in a commutative unitary ring.
\begin{definition}
Let $G$ be any group, $a,b\in G$. We write $a:b$ if $b=ac$ for some $c\in C_G(b)$ and sometimes say that $a,b$ are co-centralized in $G$. This is well-defined speech because it is a symmetric notion, since $a=bc^{-1}=c^{-1}b$, which implies that $ca=b=ac$ and hence $b:a$ with the same co-centralizer $c$.
\end{definition}
The associated relation is, then, reflexive and symmetric for arbitrary $G$, but whether it is transitive is more subtle. As a corollary of the criterion, we prove this to be the case and in fact, the relation reduces to that of commutation in $G$, which is certainly in general not transitive. To be clear, this is a distinct notion from commutative transitivity on $G$; here the relation is defined only on a proper subset of the group. However, these groups are almost commutative transitive in the sense that no elements separated by the union $(G-\gamma_2)\cup \gamma_2$ can commute unless one is central.
The first theorem in this direction shows the relationship between co-centralization and centralizers when $G$ has FL-centralizers.

Before this, we note a similar but simpler relationship for the stronger notion of \emph{co-central elements} in any group: call two elements $a,b$ co-central if they differ by a central element; let $b=az$ for some central element $z\in Z_1$ such that clearly $C_G(a)\leq C_G(b)$. If $g\in G$ commutes with $b$, then it commutes with $a=bz^{-1}$ as well such that $C_G(a)=C_G(b)$.
\begin{definition}
It is said that any group $G$ has higher free-like centralizers (or HFL-centralizers, for short) if the centralizer of any $a\in G\setminus [G,G]$ takes the form $$C_G(a)=F(b_1,\dots,b_k)\oplus Z_1$$ where the first direct factor denotes a $k\geq 1$ rank free-abelian group with generators $b_1,\dots,b_k\in G\setminus [G,G]$ of the cyclic factors in $F$.
Despite that the definition makes sense for any group, of most interest to us is when the group is finitely generated, torsion-free and nilpotent.

Note that in the case $k=1$, we can always choose as representative the generator of $F$. Express $a=b_1^iz$ for some nonzero integer $i$ and central element $z\in Z_1$. Then, $C_G(a)=C_G(b_1^i)$ since the representatives are co-central, and since $G$ is assumed torsion-free nilpotent, $C_G(a)=C_G(b_1)$.
\end{definition}
To prove the criterion we need the following well-known result on divisibility of subgroups.
\begin{proposition} (KM, [5])
A divisible subgroup of an abelian group $G$ is a direct summand of $G$.
\end{proposition} \label{divisibilityLemma}
Next, we present the quantization result for the FL-centralizer property.
\begin{theorem}
Let $G$ be a $\mathcal{T}$-group. Then $G$ has FL-centralizers if and only if for all $a, b\in G - [G,G]$, divisibility is given by $a:b \iff C_G(a)=C_G(b)$.
\end{theorem}\label{FL-CentralizerCriterion}
\begin{proof}
($\rightarrow$) Suppose first that $G$ has FL centralizers ($k=1$). If $b=ac$ for $c\in C_G(b)$, then an element $g$ in $G$ commutes with $b$ when $g=b^iz$ for some integer $i$ and $z\in Z_1$. This implies $ga=b^izbc^{-1}=bc^{-1}b^iz=ag$ trivially as all elements in the product commute with each other; $C_G(b)\leq C_G(a)$. The other inclusion holds appealing to the aforementioned symmetry. Now if $C_G(a)=C_G(b)$, then $b=a^rz=a(a^{r-1}z)$ for some integer $r$ and central element $z\in G$, where $a^{r-1}z\in C_G(b)$. Now, the case $k>1$ is nearly identical. Suppose that $b=ac$ for $c\in C_G(b)=\langle b_1\rangle\times\dots \times \langle b_k\rangle \times Z_1$, with each $b_i\in G\setminus \gamma_2$. Any element $g\in G$ that commutes with $b$ also commutes with $c$ as the centralizer $C_G(b)$ is abelian, which implies $ga=gbc^{-1}=bc^{-1}g=ag$; that is, $C_G(b)\leq C_G(a)$. In particular, then, $b\in C_G(a)$, so again by commutativity of the subgroup $C_G(a)$, we have $C_G(a)=C_G(b)$. Conversely, this equality gives $b=a_{i_1}^{k_1}\cdots a_{i_n}^{k_n}z$ for generators $a_{i_j}$ of the infinite cyclic factors of $C_G(a)$ and $z$ central in $G$. We can write this as $b=a(a^{-1}a_{i_1}^{k_1}\cdots a_{i_n}^{k_n}z)=a\zeta$ where $\zeta\in C_G(a)=C_G(b)$.
\bigskip

($\leftarrow$) Conversely, assume that any non-derived elements $a,b \in G$ are co-centralized precisely when their centralizers in $G$ coincide. Given an element $b\in G\setminus \gamma_2$, notice that for any $x\in C_G(b)$, we have $bx^{-1}:b$. It follows that $C_G(x)=C_G(b)$ for all $x\in C_G(b)$, and hence all centralizers in $G$ are abelian: for, if $x,y\in C_G(b)$, we have $$C_G(x)=C_G(b)=C_G(y).$$ Then, we can write $C_G(b)\simeq\mathbb{Z}^n$
for some integer $n\geq 1$ by the fundamental theorem of finitely generated abelian groups. Note $Z_1$ as a subgroup of $C_G(b)$ is also some power of $\mathbb{Z}$, but \emph{a priori} not a direct summand. It cannot have the same rank $n$ as $C_G(b)$ since then, it coincides with $C_G(b)$: that is, $C_G(b)=Z_1$ which means $b\in G$ is central, and then $G$ would be abelian. 

Notice that if $a^m\in Z_1$ for any $a\in G$ and $m\neq 0$, then let $k$ be the minimial positive integer such that $a\in Z_k$, and assume $k>1$. This means $aZ_{k-1}$ is a finite order element in $Z_k/Z_{k-1}$ with $(aZ_{k-1})^m=Z_{k-1}$ since $Z_{1}\subset Z_{k-1}$. The upper central factors are torsion-free, so that $a\in Z_{k-1}$, a contradiction. Hence $k=1$ and, indeed, $Z_1$ is a divisible subgroup of $G$, as well as of $C_G(b)$. By the previous lemma, $Z_1$ must be a direct summand of $C_G(b)$:
\begin{align*}
    C_G(b)=B\oplus Z_1
\end{align*}
for some subgroup $B=\langle b_1\rangle\oplus \dots\oplus \langle b_s\rangle$ of the centralizer, which must be free abelian of strictly smaller rank $s<n$. If $s=1$, clearly $b_1\notin \gamma_2$ and so this shows $G$ is an FL-centralizer group iff   $a:b$ exactly when their centralizers coincide. %It remains only to show we can select a free basis with generators from $G\setminus \gamma_2$, completing the proof. This is clear when $B$ is infinite cyclic, that is, when $G$ has FL-centralizers. The case $k>1$ could be dealt with as follows: for each $i=1,\dots,k$, write
\end{proof}
The previous result and remark immediately allow us to determine that co-centralization is indeed an equivalence relation on the subset $G-[G,G]$ of an HFL-centralizer group: for any $a,b,c\in G$, we have
\begin{align*}
    a:b \:\: \text{and} \:\: b:c \implies C_G(a)=C_G(b)=C_G(c) \implies a:c.
\end{align*}
Since $G$ has abelian centralizers, if two elements commute, their centralizers in $G$ coincide. So, it turns out this relation is equivalent, in HFL-centralizer nilpotent groups, to commutation. That is, elements satisfy $a:b$ if and only if $[a,b]=1$ when $G$ lies in this class of groups. 
\begin{corollary}
Let $G$ be a $\mathcal{T}$-group with HFL-centralizers. Then co-centralization is an equivalence relation. Moreover, two elements $a,b\in G\setminus [G,G]$ are co-centralized if and only if $[a,b]=1$ in $G$.
\end{corollary}

This equivalence is most essentially used to show that $G$ is metabelian which then aids in the classification; that is, its commutator subgroup is actually abelian. Applied to show that a $c$-nilpotent $G$ with FL-centralizers cannot be decomposed as a direct product of nilpotent groups with successive classes $1$ up to $c$. Besides this, in the next section, it is a crucial step in establishing the subgroup distortion criterion. A shorter proof of the following fact is given in a later section.
\begin{theorem} \label{ULequivalence}
Let $G$ be a $\mathcal{T}$-group with FL-centralizers of class $c\geq 2$. Then the central series of $G$ coincide.
\end{theorem}
\begin{proof}
(Step 1): If $G$ is abelian, the result trivializes. If $c=2$ and $x\in Z_{1}\setminus \gamma_{2}$, then there is $u\notin \gamma_2$ such that $G=C_G(x)=\langle u\rangle \times \gamma_2$ is abelian, contradiction. Assume then that $c\geq 3$ and for the sake of contradiction $x\in Z_{c-1}\setminus \gamma_2$. The centralizer of $x$ in $G$ takes form $\langle u\rangle \times \gamma_{c}$ for some $u\notin \gamma_2$ to which we can associate a subgroup $H(x)=\langle \gamma_{c-1},x\rangle$ of $G$. Since $[\gamma_{c-1}(G),Z_{c-1}(G)]\leq Z_0(G)=1$, $H(x)$ is abelian by construction so that $\gamma_{c-1}(G)\unlhd C_G(x)=\langle u\rangle \times \gamma_c.$
Consider the collection $$\mathcal{C}=\{H(x)|x\in Z_{c-1}\setminus \gamma_2\}$$ of all such subgroups of $G$ each satisfying $H(x)\leq C_G(x)\leq Z_{c-1}$, for $x\in Z_{c-1}\setminus\gamma_2$ using torsion-freeness of $Z_c/Z_{c-1}$. Since $\bigcup_x H(x)\subset \bigcup_x C_G(x)$, we have
\begin{align*}
    \langle H(x)|x\rangle&\leq \langle C_G(x)|x\rangle\\
    &=\langle u_1,\dots,u_n\rangle \gamma_c,
\end{align*} using the assumption of finite generation and commutation. Similarly, using the other generating set in terms of the subgroups $H(x)$,
\begin{align*}
    \langle x_1,\dots,x_n\rangle\gamma_{c-1}=\langle u_{1}^{k_1},\dots,u_{n}^{k_n}\rangle\langle z_1,\dots,z_n\rangle \gamma_{c-1}
\end{align*}
where $x_i=u_i^{k_i}z_i$ with $k_i\in\mathbb{Z}$ and $z_i\in \gamma_c$. Then
\begin{align*}
    \langle u_{1}^{k_1},\dots,u_{n}^{k_n}\rangle\gamma_{c-1}=\langle u_1,\dots,u_n\rangle \gamma_c
\end{align*}
 Note that $C_G(u_i)=C_G(x_i)$ for each $i=1,\dots,n$, since for any $\alpha=u_i^{\ell}z\in C_G(x_i)$ where $z\in \gamma_c$ central, $\alpha u_i=u_i^{\ell+1}z=u_i(u_i^{\ell}z)=u_i\alpha$ and for the converse, if $\beta \in C_G(u_i)$, then $\beta x_i=\beta u_i^{k_i}z_i=u_i^{k_i}\beta z_i=x\beta$. Then, we can repeat the construction replacing $x_i\mapsto u_i$. Indeed, this choice is valid as one can check that $H(u_i)=C_G(x_i)$ for each $i$. In fact, this held in the more general construction realizing $C_G(x_i)=C_G(u_i)\geq H(u_i)$. So
it turns out we must have
\begin{align*}
    \langle u_1,\dots,u_n\rangle \gamma_{c-1}=\langle u_1,\dots,u_n\rangle \gamma_c.
\end{align*}
Independently, suppose there exists $1\neq\beta\in Z_1$ such that $\beta=u^rz$ for $r\neq 0$ and $z\in \gamma_c$, and $u$ any element in $G\setminus \gamma_2$. Then $u^r$ is a central element of $G$ so $C_G(u)=C_G(u^r)=G$, a contradiction. This follows by the fact that in a torsion-free nilpotent group, roots are unique when they exist, and so
\begin{align*}
 \alpha^nb^m&=b^m\alpha^n, \: (b^{-m}\alpha b^m)^n=\alpha^n, \:
 b^{-m}\alpha b^m=\alpha, \:
 b^m=(\alpha b\alpha^{-1})^m, \: \alpha=\alpha^{b}
\end{align*} for $n,m\neq 0$. Hence $r=0$, and we conclude $Z_1=\gamma_c$.

If $\langle u_1,\dots,u_n\rangle \cap \gamma_{c-1}=1$, then $\langle u_1,\dots,u_n\rangle\cap Z_1=1$ and a standard fact concludes that $\gamma_{c-1}=Z_1$, which cannot be. So, there is $1\neq x=u_1^{k_1}\dots u_n^{k_n}\in \gamma_{c-1}\setminus Z_1$ and recall we still have $\gamma_{c-1}\subseteq C_G(u_k)=\langle u_k\rangle \times Z_1$ for each $k$. For any $i\neq j$, $x\in C_G(u_i,u_j)$, so that $x=u_i^{m_i}\nu_i=u_j^{m_j}\nu_j$ for integers $m_i,m_j\neq 0$ and $\nu_i,\nu_j\in Z_1$. Then $u_i^{m_i}\in C_G(u_j)$ and the same reasoning gives $C_G(u_i)=C_G(u_j)$. Since $1\leq i,j\leq n$ are arbitrary, we have $A=C_G(u_i)$ for all $i=1,\dots,n$ and hence 
\begin{align*}
    Z_{c-1}=A\cup \gamma_2.
\end{align*} If $Z_{c-1}=\gamma_2$, we are done, so assume $Z_{c-1}=A$, an abelian subgroup of $G$.
Since $Z_{c-j}$ is abelian and $Z_{c-j}/Z_1\simeq \mathbb{Z}$ for each $j<c-1$, then $Z_{c-j}=\langle u^{r_j}\rangle \oplus Z_1$ for some integers $1=r_1\leq r_2\leq \dots\leq r_{c-2}$ and $r_{c-1}=0$. But 
\begin{align*}
    Z_2/Z_1=\langle u^{r_{c-2}}\rangle \oplus Z_1/Z_1\simeq \langle u^{r_{c-2}}\rangle/\langle u^{r_{c-2}}\rangle \cap Z_1=\langle u^{r_{c-2}}\rangle.
\end{align*} The second upper central factor is
\begin{align*}
    Z_3/Z_2\simeq \frac{Z_3/Z_1}{Z_2/Z_1}\simeq \langle u^{r_{c-3}}\rangle/\langle u^{r_{c-2}}\rangle.
\end{align*}
This is isomorphic to $\mathbb{Z}_q$ for some $q\neq 0$ or is trivial corresponding to whether $r_{c-3}<r_{c-2}$ or $r_{c-3}=r_{c-2}$. Neither is acceptable, either contradicting torsion-freeness or the nilpotency class of $G$.

\bigskip

(Step 2): For any $2\leq i<c$, suppose there exists $x\in Z_{c-i+1}\setminus \gamma_i$ and by induction that $Z_{c-j+1}=\gamma_j$ for all $j<i$. Note $x\in \gamma_{i-1}$ is a product of $(i-1)$-commutators, not all of which can be $i$ commutators- have a $2$-commutator "nested" in at least one component. That is, 
\begin{align*}
    x=\Pi_{i=1}^m[a_i,b_i],
\end{align*}
where there must exist $i=1,\dots,m$ such that $a_i\notin \gamma_2$ and $b_i\notin \gamma_{i-2}$. Taking right weighted commutators, we have $b_i\in \gamma_{i-2}$ for each $i$, making the second condition superfluous. This commutator exists otherwise each commutator in the product would be $i$-step. Let $A_x\subseteq G\setminus \gamma_2$ consist of the first components of the pairs $(a_i,b_i)$, the commutators of the elements of which product to $x$, as shown above, and $a_i\notin \gamma_2$. Note, for a given $x\in Z_{c-i+1}\setminus \gamma_i$, the set $A_x$ is finite. Associated to each $a\in A_x$ is $u_{x,a}\notin \gamma_2$ such that $C_G(a)=\langle u_{x,a}\rangle \oplus Z_1$. Consider the subgroup 
\begin{align*}
    H&=\langle C_G(a): a\in A_x, x\in Z_{c-i+1}\setminus \gamma_2\rangle\\
    &=\langle  u_{x,a}:  a\in A_x, x\in Z_{c-i+1}\setminus \gamma_2 \rangle Z_1\\
    &=\langle  u_{x_i,a^i_j}: x_i\in Z_{c-i+1}, a^i_{1},\dots,a^i_{k_i}\in A_{x_i}\rangle Z_1
\end{align*}
of $G$. Note that the product is not necessarily semi-direct. One may try to use the argument that showed $C_G(u_i)=C_G(u_j)$ for all $i\neq j$ in step 1, but then we assumed that $x\notin Z_1$. %The intersection of the two subgroups in the product is trivial, due to the special form of centralizers in $G$. Explicitly, suppose $y=u_1^{k_1}\dots u_n^{k_n}\in Z_1$ with $k_i\neq 0$ is a nontrivial element of the intersection. For ease, adopt the notation $h_{ij}=u_i^{k_i}\dots u_j^{k_j}$. Note that $y=u_1^{k_1}h_{2n}\in C_G(u_1)$ means $h_{2n}\in C_G(u_1)$, so that
%\begin{align*}
    %h_{2n}=u_1^{r_1}z_1=u_2^{k_2}\dots u_n^{k_n}\\
    %u_1^{-r_1}h_{2n}=z_1\in Z_1
%\end{align*} where $z_1\in Z_1$ and $r_1$ is an integer. If $h_{2n}\in Z_1$, then $u_1^{k_n}\in Z_1$, a contradiction. 

Next, assume there exists $b\in G\setminus (\gamma_2\cup H)$ with centralizer $C_G(b)=\langle v\rangle \oplus Z_1$ and write $b=v^iz$ with $i\neq 0$ and $z\in Z_1$, implying $v\notin H$. In particular,
\begin{align*}
    v\in \bigcap_{x\in Z_{c-i+1}\setminus \gamma_i}(G\setminus A_x).
\end{align*}
Next, there are two cases to consider. Suppose $v$ appears in the first component of \emph{some} pair of elements whose commutator is a non-trivial factor in the decomposition of an element $y\in Z_{c-i+1}\setminus \gamma_i$ into $(i-1)$-commutators. Then $v\notin H$ implies $v\in \gamma_2$, since the second component in any $[\alpha,\beta]\in \gamma_{i-1}$ will always be $(i-2)$-commutator, as remarked above, a contradiction. In the second case, it follows that $v\in C_G(\gamma_{i-2})$, which means $Z_{r}\leq \gamma_{i-2}\subseteq C_G(b)$ for all $r\leq c-i+3$, by the induction hypothesis. Because $c-i+3\geq 3$, the same argument given in step 1 shows, for example, that $Z_3/Z_2$ is a torsion group, a contradiction.

Hence, we have $G=\gamma_2\cup H$ implying $G=H$ since the derived subgroup is improper when $c\geq 2$. But this cannot be, for then let $\{g_1Z_{c-1},\dots,g_{s}Z_{c-1}\}$ be a free basis of the free abelian group 
\begin{align*}
    Z_c/Z_{c-1}=\langle g_1Z_{c-1}\rangle \oplus \dots \oplus \langle g_sZ_{c-1}\rangle
\end{align*} For a given $i=1,\dots,s$, we can express $g_i=\sum_j m_ju_j$ where for each $j$, we have $u_j=u_{x_p,a_q^p}$ for some integers $p,q,m_j$. Note there must be a nontrivial element $\sum_{i} r_iu_i\in Z_{c-1}$. Write each element $r_iu_i+Z_{c-1}$ in the free basis as $\sum_{k} r_{i,k}'g_k+Z_{c-1}$ for possibly zero integers $r_{i,k}'$. Then,
\begin{align*}
    Z_{c-1}&=(\sum_i r_iu_i)+Z_{c-1}=\sum_i r_i(u_i+Z_{c-1})\\
    &=\sum_i\sum_k r_{i,k}'(g_k+Z_{c-1})
\end{align*}which implies $r_{i,k}'=0$ for all $i,k$. Hence, for any $i$, we have $r_iu_i\in Z_{c-1}$ and by torsion-freeness $u_i\in Z_{c-1}$, a contradiction. Thus, concluding the inductive step.
\end{proof}
Before proceeding, we state some useful terminology that may not be convention.
\begin{definition}
Let $G$ be any nilpotent group. If the upper central series of any upper central term $Z_k$ consists of the remaining terms in that of $G$,
\begin{align*}
    Z_k\geq Z_{k-1}\geq \dots\geq Z_1\geq Z_0=1,
\end{align*}
then we will say $G$ is said to be \emph{(strongly) tight with respect to nilpotency class}. When this occurs, the higher centers of $G$ satisfy many special properties. For instance, the upper central term $Z_k$ of $G$ has precisely class $k$ and center coinciding with that of $G$, for all $k=1,\dots,c$. In this case, we will say $G$ is \emph{tight with respect to its class}. These facts follow directly from the definition as the UCS of $G$ has length $k$.
\end{definition}
It is known that for any integer $c\geq 2$, there exists a group that is tight with respect to its nilpotency class. We prove this here for reference and because the construction in its proof is indirectly used in the proceeding corollary.
\begin{lemma}
For any $c\geq 2$, there exists a $c$-nilpotent group $G$ that is tight with respect to its nilpotency class.
\end{lemma}
\begin{proof}
Select any nilpotent groups $K_1,\dots,K_c$ with each $K_i$ having class $i$ and let $G=K_1\times \dots\times K_c$ be the external direct product of them. Then, as a direct product of nilpotent groups, $G$ is nilpotent. Moreover, its nilpotency class is $G$ using the fact that
\begin{align*}
    Z_k(G)=Z_k(K_1)\times \dots\times Z_k(K_c)
\end{align*}
so $c$ is the smallest integer such that $Z_c(G)=K_1\times \dots\times K_c=G$. More generally, for any $i=1,\dots,c$, we have
\begin{align*}
    Z_i(G)=K_1\times \dots \times K_i\times Z_{i}(K_{i+1})\times \dots\times Z_i(K_{c}).
\end{align*}
Examining the LCS of $Z_i$ and its component groups, the class of $Z_i(G)$ cannot drop below $i$, the class of $K_i$. Also, the class of $Z_i$ cannot exceed $i$ because the components have at most class $i$. So, either by again taking repeated commutators or quotienting by their centers and taking higher centers of these factor groups shows $\text{Class}(Z_i(G))=i$ for each $i$.
\end{proof}
\begin{corollary}
 Let $G$ be a $c$-nilpotent group with FL-centralizers, $c\geq 3$. Then the group $G$ is indecomposable as a product of groups
 \begin{align*}
     G_1\times \dots \times G_c,
 \end{align*}
 where each $G_i$ is nilpotent of class $i$.
\begin{proof}
By the proof of the previous proposition, it suffices that $G$ is not tight with respect to its nilpotency class, $c$. For this, we produce at least one upper central term $Z_k$ with class less than $k$. Select any term $\gamma_j$ in the lower half of the LCS of $G$, so $c>j>c/2$ if $c$ is an even integer and $c>j\geq \frac{c-1}{2}$ if $c$ is odd. Note $[\gamma_j,\gamma_j]\leq \gamma_{2j}=\{1\}$. This together with Theorem \ref{ULequivalence} gives that $Z_{c-j+1}$ is abelian, while $Z_{c-j+1}\neq Z_1$.
\end{proof}
\end{corollary}

\begin{thebibliography}{}
{\Large
\bibitem{Roots of Nilpotent Groups}
  \textit{Sze, S., Kahrobaei, D., Dambreville, R., and Dupas, M. (2011). Finding N-th Roots in Nilpotent Groups and Applications to Cryptology. International journal of pure and applied mathematics, 70, 571.}\\
\bibitem{Osin}
  \textit{Osin, D.V. Subgroup Distortions in Nilpotent Groups, Communications in
 Algebra, Vol. 29, Issue 12, 5439-5463. (2010)}\\
\bibitem{Robinson} 
   \textit{ D.J.S. Robinson. A Course in the Theory of Groups, Graduate Texts in Mathematics; Springer, Berlin, GTM 80. (1996)}\\
\bibitem{Warfield}   
   \textit{Warfield, R. B. Sections 1-3. In Nilpotent groups. essay, Berlin, Heidelberg, New York: Springer. (1976)}\\
\bibitem{KM}
    \textit{ M. I. Kargapolov and J. I. Merzlejakov. Fundamentals of the Theory of
 Groups, Graduate Texts in Mathematics; 62, Springer Verlag New York
 Inc. (1979)}\\
\bibitem{Automorphisms}
    \textit{ Escayola, Maximiliano and Rivas, Cristóbal. (2023). On the critical regularity of nilpotent groups acting on the interval: the metabelian case. 10.48550/arXiv.2305.00342. }
}

\end{thebibliography}
\end{document}